\documentclass[12pt]{amsart}
\usepackage {amsxtra,amssymb,amsmath,amsfonts,amscd}

\theoremstyle{plain}
\newtheorem{theorem}{Theorem}[subsection]
\newtheorem{corollary}[theorem]{Corollary}
\newtheorem{lemma}[theorem]{Lemma}
\newtheorem{proposition}[theorem]{Proposition}
\newtheorem{definition}[theorem]{Definition}
\newtheorem{remark}[theorem]{Remark}
\newtheorem{note}[theorem]{Note}

\numberwithin{equation}{subsection}

\newcommand{\fieldk}{\ensuremath{\Bbbk} }
\newcommand{\tensor}{\ensuremath{\bigotimes} }
\begin{document}
\title[Differential Operators]{Differential Operators on Azumaya algebra and 
          Heisenberg algebra}
\author{Uma N. Iyer}
\address{Mehta Research Institute of Mathematics
	 and Mathematical Physics,
          Chhatnag Road, Jhusi\\
         Allahabad - 211019, U.P., India}
\email{uiyer@mri.ernet.in}
\date{\today }
\maketitle
\section{Introduction}\label{S:intro}
Let \fieldk be a field  and $A$ an associative \fieldk -algebra.  Let 
\tensor mean $\tensor _{\fieldk }$ and $A^e := A \tensor A^o$, where
$A^o$ denotes the opposite algebra of $A$.  Any left $A^e$-module $M$ is the
same as an $A$-bimodule given by $a\cdot m\cdot b = (a\otimes b^o)\cdot m$
for $m\in M$ and $a,b \in A$.
In \cite{LR}, V. Lunts and A. Rosenberg give a definition for the ring of
differential operators on a left $A$-module $L$ (denoted by
$D_{\fieldk }(_AL)$) and the $A$-bimodule of differential operators on $L$ of
order less than or equal to $m$ denoted by $D_{\fieldk}^m(_AL)$.
When $L= A$, we denote the ring and module respectively by
$D_{\fieldk}(A)$ and $D^m_{\fieldk}(A)$.  These definitions are recalled 
in the section of preliminaries in this paper.

In this paper we compute the ring of differential operators for some
noncommutative rings, namely the Axumaya algebras and the Heisenberg algebras.
The initial interest was to compute these rings for matrix algebras and
the Weyl algebras, which were easily generalized.

We consider Azumaya
algebras over a noetherian ring.  
We conclude that the ring of differential operators are generated as modules
by the ring of differential operators on their centre and 
homomorphisms given by multiplication by elements of the bigger ring
(called inner homomorphisms). 
That is, if 
$R$ is the centre of $A$ (where $A$ is an Azumaya algebra over $R$), 
we show that $D_{\fieldk}(R)$ can be embedded into
$D_{\fieldk}(A)$ and that 
$D_{\fieldk}(A) = (A\tensor _R A^o).D_{\fieldk}(R).(A\tensor _R A^o)$, 
that is, $D_{\fieldk}(A)$ is generated
as an $A$-bimodule by $D_{\fieldk}(R)$ (Theorem \ref{T:DA=ADRA}).

In the case of Heisenberg algebras of zero characteristic, 
we need two copies of
differential operators on the centre to generate all the differential 
operators (Theorem \ref{T:Hnchar0}).  The non zero characteristic
follows from the study on Azumaya algebras, because in this case the Heisenberg
algebra is Azumaya over its centre (Theorem 2 \cite{R}).

In particular, our work on these general rings show that
\begin{enumerate}
\item	If  $R$ is a \fieldk -algebra, we show that
	$D_{\fieldk }(M_n(R)) = M_{n^2}(D_{\fieldk }(R))$, where
	$M_n(R)$ denotes the algebra of $n \times n$ matrices over $R$
	(Corollary \ref{C:mat}).
\item	If $A_n$ denotes the $n$-th Weyl algebra over a field of characteristic
	0, then $D_{\fieldk }(A_n) = A_{2n}$ (Corollary \ref{C:weyl}).
\end{enumerate}
In the case of Azumaya algebras, we show
that there is a one-to-one correspondence between ideals of $D_{\fieldk }(A)$
and $D_{\fieldk }(R)$ (section \ref{S:ideals}).
If $H_n$ denotes the $n$th-
Heisenberg algebra, we show that $D_{\fieldk}(H_n)$ is simple 
(Theorem \ref{T:DHnsimple} and corollary \ref{C:Hncharp}).

We give some definitions and
prove some elementary results in the section of preliminaries.  These 
results will be used later, and are interesting in their own right.
This will be followed by a section on the differential operators on the Azumaya
algebras.  Here, we first show that if $A$ is an Azumaya algebra over $R$,
then $D^m_{\fieldk}(A) = D^m_{\fieldk}(_RA)$ (Theorem \ref{T:DA=DRA})
for each $m\geq 0$. Then we show that
$D_{\fieldk}(R)$ embeds as an $R$-bimodule
in $D_{\fieldk}(_RA)$ (respecting the filtration
given by the order of differential operators), and along with the inner 
differential operators, generate the entire ring $D_{\fieldk}(_RA)$
(Theorem \ref{T:DA=ADRA}).

The last section covers the Heisenberg algebras.  We consider the two cases
of zero characteristic and non zero characteristic separately.  \\
\textbf{Notations :}
\begin{enumerate}
\item	For any $\Bbbk$-algebra
 	$S$, and $s \in S$, denote by $\lambda _s \in 
	\mathit{Hom}_{\Bbbk}(S,S)$\\
	(respectively $\rho _s$) to be the homomorphism 
	given by left-multiplication
 	(respectively, right-multiplication) by $s$. 
\item	For any $\varphi \in \mathit{Hom}_{\Bbbk}(S,S)$, let
	$[\varphi , s] := \varphi \cdot s - s \cdot \varphi$.
\item 	For $r,s \in S$, let $[r,s] := rs -sr$.
\end{enumerate}
\section{Preliminaries}
We recall from \cite{LR} some definitions.
\begin{definition}\label{D:diff} \hfill
\begin{enumerate}
\item	For an $A^e$ module $M$, its centre is the \fieldk -submodule
	\[
	\mathcal{Z}(M) := \{ z\in M | a\cdot z = z \cdot a \textit{ for }
	a \in A \}.
	\]
\item	Define the $i$-th differential part of $M$,
	$\mathcal{Z}_iM$ by induction as follows:
	\begin{align*}
	\mathcal{Z}_0M &:= A^e \mathcal{Z}(M), \textit{ and}\\
	\mathcal{Z}_iM \diagup \mathcal{Z}_{i-1}M &:= A^e \mathcal{Z}\left( 
	M\diagup
	\mathcal{Z}_{i-1}M \right) \textit{ for } i \geq 1.
	\end{align*}
\item	The differential part of an $A^e$ module $M$ is
	$M_{\textit{diff}} := \cup_{i\geq 0}\mathcal{Z}_iM$.
\end{enumerate}
\end{definition}
For $L$ a left $A$-module, the \fieldk -vector space 
$\mathit{Hom}_{\fieldk} (L,L)$ has an $A^e$-module
structure given by $c\cdot \varphi \cdot a (b) = c\varphi (ab)$ for
$a,b,c \in A$ and $\varphi \in \mathit{Hom}_{\fieldk} (L,L)$.
The \textit{differential part} of $\mathit{Hom}_{\fieldk}(L,L)$ is
the \textit{algebra of \fieldk -linear differential operators on $L$}, and
denoted by $D_{\fieldk}(_AL)$.
The $A^e$ module of differential operators of order $\leq m$ on $A$ is
$\mathcal{Z}_m \mathit{Hom}_{\fieldk}(L,L)$ and is denoted by
$D^m_{\fieldk}(_AL)$.  This definition generalizes the definition of
differential operators on a commutative ring as given by Grothendieck 
(\cite{G2}).
We denote $D_{\fieldk}(_A A)$ simply by $D_{\fieldk}(A)$.

We state and prove some preliminary results.
\begin{proposition}
For any ring $A$, the ring $D^0_{\fieldk}(A)$ consists of left and right 
multiplications by elements of $A$.
\end{proposition}
\begin{proof}
The central elements of the $A^e$-module $\mathit{Hom}_{\fieldk}(A,A)$ 
are homomorphisms given by right multiplication by elements of $A$.
Hence the result.
\end{proof}
\begin{corollary}\label{C:D0}
There is a surjection 
\begin{align*}
A\tensor _{\mathcal{Z}(A)} A^o &\to D^0_{\fieldk}(A) \textit{ given by}\\
a\otimes b^o &\mapsto [ c \mapsto acb],
\end{align*}
where $\mathcal{Z}(A)$ is the centre of the ring $A$.
\end{corollary}
\begin{proposition}\label{P:bigtocentre}
Let $R \subset S$ be two \fieldk -algebras, and $M$ be an $S^e$-module
(hence an $R^e$-module).  If $R\subset \mathcal{Z}(S)$, the centre of $S$,
then the $i$-th differential part of $M$ considered as an 
$S$-bimodule is contained in the
$i$-th differential part of $M$ considered as an $R$-bimodule.
\end{proposition}
\begin{proof}
For any $S$-bimodule $N$
	we have,
\[
S\cdot \mathcal{Z}_S (N) \subset \mathcal{Z}_R(N),
\]
where $\mathcal{Z}_S(N)$ denotes the $S$-centre of $N$ (analogously defined
for $R$).  Hence the proposition.
\end{proof}
\begin{corollary}\label{C:bigtocentre}
Let $R\subset S$ be two \fieldk -algebras.  
If $R\subset \mathcal{Z}(S)$, then 
$D^m_{\fieldk }(S) \subset D^m_{\fieldk }(_RS)$ for $m\geq 0$.
\end{corollary}
\begin{remark}
The corollary above is not true if $R\nsubseteq \mathcal{Z}(S)$.  Consider for 
example,
\[
R= \Bbbk [x] \subset S = \Bbbk <x,y> \diagup <[y,x] = y>.
\]
Note that $R$ is commutative.  Hence, $\varphi \in D_{\Bbbk}(_RS)$
satisfies
\[
[[\cdots [ \varphi ,r_1],r_s],\cdots r_n] = 0,
\]
for some $n\geq 0$.
We know that $\lambda _y$, the homomorphism given by left multiplication by
$y$ is in $D_\Bbbk (S)$.  
But $[\lambda _y ,x] =  \lambda _y$.  Hence
$\lambda _y \notin D_{\Bbbk }(_RS)$.
\end{remark}
Let $L$ be a free, left-$R$-module, where $R$ is a commutative 
\fieldk -algebra.  Fix a basis $\{ l_1, l_2,\cdots , l_n \}$ of $L$ over $R$.
Any $\Phi \in \mathit{Hom}_{\fieldk }(A,A)$ can be written as 
\begin{equation}\label{E:homfree}
\Phi =
	\begin{matrix}
		& \begin{array}{cccc} 
			R.l_1 & R.l_2 & \cdots & R.l_n 
	     	   \end{array}\\
	     \begin{array}{llll}
		R.l_1 \\ R.l_2 \\ \vdots \\ R.l_n 
	     \end{array}
		&
		\begin{pmatrix}
		\varphi _{1,1} & \varphi _{1,2} & \dots & \varphi_{1,n}\\
		\varphi _{2,1} & \varphi _{2,2} & \dots & \varphi_{2,n}\\
		\vdots 	       & \vdots         & \ddots & \vdots \\
		\varphi _{n,1} & \varphi _{n,2} & \dots & \varphi _{n,n}
		\end{pmatrix}
	\end{matrix}
\end{equation}
where $\varphi _{i,j} \in \mathit{Hom}_{\fieldk}(R,R)$.
\begin{proposition}\label{P:free}
Referring to the equation \ref{E:homfree}, $\Phi \in D^m_{\fieldk}(_RL)$
if and only if $\varphi _{i,j} \in D^m_{\fieldk}(R)$.
\end{proposition}
\begin{proof}
Since $R$ is commutative,\\
$\Phi \in D^m_{\fieldk}(_RA)$ (respectively, 
$\varphi \in D^m_{\fieldk}(R)$), if and only if \\
$[\cdots [[\Phi , r_0],r_1],\cdots ,r_m] =0$ 
(respectively,
$[\cdots [[\varphi , r_0],r_1],\cdots ,r_m] =0$), for 
$r_0,r_1,\cdots , r_m \in R$.
The proposition follows immediately once we notice that 
if $\Phi$ is given by a matrix $(\varphi _{i,j})$, then
$[\Phi ,r]$ is given by the matrix $([\varphi _{i,j},r])$ for $r\in R$.
\end{proof}
\section{Differential operators on Azumaya algebras}\label{S:azumaya}
Let $R$ be a commutative, Noetherian \fieldk -algebra.  
Let $A$ be an Azumaya algebra over $R$ (see \cite{DI} for a complete study); 
i.e., $A$ is an $R$-algebra which is
finitely generated,
projective, and faithful as an $R$-module, such that $R\cdot 1 = 
\mathcal{Z}(A)$
and the map
\begin{align*}
A\tensor _R A^o &\to \mathit{Hom}_R(A,A), \\
a\otimes b^o &\mapsto [c \mapsto acb]
\end{align*}
is an isomorphism.
Examples are matrix algebras over $R$.
Some immediate remarks follow:
\begin{remark}\label{R:AzD0}
$D^0_{\fieldk}(A) = \mathit{Hom}_R(A,A) \cong A \tensor _R A^o$, and hence
\\
$D^0_{\fieldk}(A) = D^0_{\fieldk}(_RA)$.
Indeed, referring to the corollary \ref{C:D0}, there is a surjection
\begin{equation}\label{E:AzD0}
 	A \tensor _RA^o \to D^0_{\fieldk}(A)
\end{equation}
On the other hand, since $A$ is an Azumaya algebra,\\
 $A \tensor _R A^o \cong
\mathit{Hom}_R(A,A)$ given by the map $a\otimes b^o (c) = acb$. Thus,
$A\tensor _RA^o$ injects into $\mathit{Hom}_{\fieldk}(A,A)$.  Hence, the 
surjection \ref{E:AzD0} onto $D^0_{\fieldk}(A)$ is an isomorphism.
By definition, $D^0_{\fieldk}(_RA) = \mathit{Hom}_R(A,A)$. 
\end{remark}
\begin{remark}
By corollary \ref{C:bigtocentre}, for each $m\geq 0$, we have a map of
$R$-bimodules, namely
\begin{equation}\label{E:DAinclDRA}
\iota _A : D^m _{\fieldk}(A) \hookrightarrow D^m_{\fieldk}(_RA).
\end{equation}
\end{remark}
\subsection{Proof of $D_{\fieldk}(_RA) = D_{\fieldk}(A)$}
\begin{theorem}\label{T:DA=DRA}
The inclusion of \ref{E:DAinclDRA} is an isomorphism.  That is,
for each $m\geq 0$, we have
$D^m_{\fieldk}(A) = D^m_{\fieldk}(_RA)$
as $R$-bimodules.
\end{theorem}
\begin{proof}
We first prove the theorem in the case when $A$ is free over $R$ with basis
$\{a_1 ,a_2, \cdots , a_n \}$.  By proposition \ref{P:free}, any $\Phi \in 
D^m_{\fieldk}(_RA)$ if and only if all the 
$\varphi _{i,j} \in D^m_{\fieldk}(R)$.  It remains to show that if
all the $\varphi _{i,j} \in D^m_{\fieldk}(R)$, then 
$\Phi \in D^m_{\fieldk}(A)$.  
Let $a_i \cdot a_j = \sum _k r_{i,j}^k a_k$. For any $\varphi \in 
\mathit{Hom}_{\fieldk}(R,R)$, define
$\tilde{\varphi} \in \mathit{Hom}_{\fieldk}(A,A)$ as
$\tilde{\varphi} (r a_i) = \varphi (r) a_i$.
For each $1 \leq l,k \leq n$, define $\wp ^{l,k} 
\in \mathit{Hom}_R(A,A)$ (and hence in $D^0_{\fieldk}(A)$) given by
$\wp ^{l,k} (a_i) = \delta _{i,k} a_l$.
Then, we have $\Phi = \sum _{i,j} \wp ^{i,1} \widetilde{\varphi _{i,j}}
\wp ^{1,j}$.
Thus, it remains to show that if $\varphi \in D_{\fieldk}^m (R)$, then
$\tilde{\varphi} \in D_{\fieldk}^m (A)$.
Using induction  on $m$ and the following identity,
\[
[\tilde{\varphi} , a_j] = \sum _{i,k} \widetilde{[\varphi , r_{j,i}^k]}
			\wp ^{k,i}
\]
we conclude the theorem in the case when $A$ is free as an $R$-module.

In the case when $A$ is not free as an $R$-module, we consider the localization
of $A$ with respect to a prime ideal $P$ of $R$.
By Lemma 5.1, pg61 of \cite{DI}, $A_P := R_P \tensor _R A$ is an Azumaya
$R_P$-algebra. Consider the injective (by flatness of $R_P$ as 
an $R$-module) map
\begin{equation}\label{E:local}
id \otimes \iota _A : R_P \tensor _R D^m _{\fieldk}(A) \to
                       R_P \tensor _R D^m _{\fieldk}(_RA).
\end{equation}
By Proposition 16.8.6 of \cite{G2}, $R_P \tensor _R D^m _{\fieldk}(_RA)
\cong D^m_{\fieldk}(_{R_P}A_P)$ which by our discussion on the free Azumaya 
case is isomorphic to $D^m_{\fieldk}(A_P)$.
Thus, it is sufficient to show that the inclusion of equation \ref{E:local}
$id \otimes \iota _A :R_P \tensor _R D^m _{\fieldk}(A) \to
	D^m _{\fieldk}(A_P)$ is surjective.  The following lemma proves this
which completes the theorem.
\end{proof}
\begin{lemma}
For $m \geqslant 0$,
the map 
\[
\mathit{id} \otimes \imath_A : R_P \bigotimes _R
            	 D_\Bbbk^m(A) \longrightarrow 
		 D^m _{\Bbbk }(A_P)
\]
            	 is surjective. 
\end{lemma}
\begin{proof}
We prove both the statement by induction on
$m$.  Let 
\[ x = \sum _i (a_i / s_i ) \otimes (b _i / t _i)^o 
\in (A_P \bigotimes _{R_P} A_P^o) \cong  D_\Bbbk ^0 (A_P)
\]
be given.
There is an $s$ in $R \setminus P$, such that
$ sx = \sum _i a'_i  \otimes (b' _i )^o \in D_\Bbbk ^0
(A)$.
Thus, 
$  (1/s)\otimes sx \in R_P \bigotimes _R D _\Bbbk ^0 (A)$
is mapped to $x$ under $\mathit{id} \otimes \imath _A$.  
So, the result is proved for $m = 0$.\\

Assuming that the proposition is proved for $m$
(which implies that $D^m_{\fieldk}(A) = D^m_{\fieldk}(_RA)$),
let 
\[
d \in R_P \bigotimes _R D_\Bbbk ^{m+1}({}_R A) =
D_\Bbbk ^{m+1} ({}_{R_P} A_P)=
D_\Bbbk ^{m+1}(A_P),
\]
be such that 
$(a/s)\cdot d - d \cdot (a/s) \in D_\Bbbk ^{m} (A_P)$ 
for every $(a/s) \in A_P$.  
It is enough to show that $d$ is in the image of
$(\mathit{id} \otimes \imath _A)$. 
Note,
$s d \in D_\Bbbk ^{m+1} ({} _R A)$ for some
$s \in R \setminus P$.  Let 
$\{ a_1 , a_2 , \cdots , a_n\}$ be a finite set of generators of $A$
as an $R$-module. 
$(a_i/1)\cdot (sd) - (sd) \cdot (a_i/1) \in 
D_\Bbbk ^{m}({}_{R_P}A_P)$ .
By induction hypothesis, for each $i$,  $1 \leqslant i \leqslant n$,
there exists a $t_i$ in $R \setminus P$, such that
$t_i \cdot [a_i \cdot (sd) - (sd) \cdot a_i] \in D_\Bbbk ^m
(A)$.
Let, $ t = t_1 t_2 \cdots t_n$.
Then, 
$[a_i \cdot (tsd) - (tsd) \cdot a_i] \in D_\Bbbk ^m (A)$, for
all
$i$, $1 \leqslant i \leqslant n$.
Let $ts = T \in R$.  For, $r \in R$ and $a_i$ a generator, consider
$(ra_i) \cdot (Td) - (Td)\cdot (ra_i) = 
r[a_i \cdot (Td) - (Td)\cdot a_i] + 
[r\cdot (Td) -(Td)\cdot r]\cdot a_i$. 
Now, 
$r[a_i \cdot (Td) - (Td)\cdot a_i] \in D_\Bbbk ^m (A)$.
Since  $sd \in D_\Bbbk ^{m+1} ({}_R A)$,
$Td \in D_\Bbbk ^{m+1} ({}_R A)$,
which implies 
$[r\cdot (Td) -(Td)\cdot r] \in D_\Bbbk ^m ({}_R A)$. 
But,  by induction hypothesis, 
$D_\Bbbk ^m ({}_R A) = D_\Bbbk ^m (A)$. Hence,
$[r\cdot (Td) -(Td)\cdot r]\cdot a_i \in D_\Bbbk ^m (A)$.
Hence, for any $a \in A$, 
$a\cdot (Td) -(Td)\cdot a \in D_\Bbbk ^m (A)$.
Thus, $(Td) \in D _\Bbbk ^{m+1}(A)$.  
Hence, $d \in R_P \bigotimes _R D _\Bbbk ^{m+1}(A)$.  
This proves the lemma.
\end{proof}
\begin{corollary}\label{C:mat}
Let $M_n(R)$ denote the ring of matrices over $R$ where
$R$ is a commutative 
\fieldk -algebra. Then,
\begin{align*}
D^m_{\fieldk}(M_n(R)) &= M_{n^2}(D_{\fieldk}^m(R)),\textit{ and hence }\\
D_{\fieldk}(M_n(R)) &= M_{n^2}(D_{\fieldk}(R)).
\end{align*}
\end{corollary}
\begin{proof}
 The theorem above shows that
$D^m_{\fieldk}(M_n(R)) = D^m_{\fieldk}(_RM_n(R))$.
The ring $M_n(R)$ is free as a left $R$-module. By Proposition \ref{P:free}
the corollary is proved.
\end{proof}
\begin{remark}
In \cite{I} we have proved a more general statement.
If $R$ and $S$ are two \fieldk -algebras such that $S$ is finite dimensional
as a \fieldk -vector space, then $D_{\fieldk}(R \tensor S) = 
D_{\fieldk}(R) \tensor D_{\fieldk}(S)$.
\end{remark}
\subsection{$D_{\fieldk}(A)$ is generated by $D_{\fieldk}(R)$ and 
		inner homomorphisms.}
Here we embed $D^m_{\fieldk}(R)$ (as $R$-bimodules) into
$D^m_{\fieldk}(_RA)$ for each $m\geq 0$.

By Lemma 3.1 of \cite{DI}, $R$ is an $R$-direct summand of $A$;
that is, $A \cong R \oplus B$ as left $R$-modules.
Fix one such decomposition.
Since $A$ is projective as a left $R$-module, $B$ is also a projective as  a
left $R$-module. 
By assumption, $R$ is a Noetherian ring.  Hence $B$ is a finitely generated
$R$-module.  By the Dual Basis Lemma (lemma 1.3 of \cite{DI}), we choose a 
a collection $\{ b_i , f_i \}_{1 \leq i \leq n}$, where $b_i \in B$ and
$f_i \in \mathit{Hom}_R(A,R)$ (we can consider $f_i$ to be elements of
$\mathit{Hom}_R(A,A)$ by the natural inclusion of $R$ into $A$)
such that $b = \sum _i f_i(b) b_i$ for $b\in B$.
Let $f_0 \in \mathit{Hom}_R(A,A)$ be the projection of $A$ onto $R$, and
$b_0 = 1$.  Extend $f_i$ for $i\geq 1$ to $A$ by defining $f_i (r)=0$ (we 
denote the extension also by $f_i$).
 Then the collection $\{b_i, f_i \}_{0\leq i \leq n}$ is
 a dual basis of $A$.
\begin{remark}
By definition, for $0\leq i \leq n$, the homomorphisms $f_i$ are
differential operators of order 0.  That is, they are inner homomorphisms.
\end{remark}
We describe a way to extend elements of $\mathit{Hom}_{\fieldk}(R,R)$ to that 
of $\mathit{Hom}_{\fieldk}(A,A)$.
\begin{definition}\label{D:ext}
For $\varphi \in \mathit{Hom}_{\fieldk}(R,R)$, define $\bar{\varphi} \in
\mathit{Hom}_{\fieldk}(A,A)$ as 
\[
\overline{\varphi} = \sum _{i=0}^n \rho _{b_i} \varphi f_i,
\]
where $\rho _{b_i}$ is the homomorphism given by right multiplication by
$b_i$.
\end{definition}
Since $a = \sum_{i \geq 0} f_i(a) b_i$, we have,
$\overline{id} = id$ where $id$ denotes the identity homomorphism in
the respective rings.
An immediate consequence is the following lemma.
\begin{lemma}
If $\varphi \in D^m_{\fieldk}(R)$, then $\overline{\varphi } 
\in D^m_{\fieldk}(A)$.
\end{lemma}
\begin{proof}
It is clear to see that 
$\overline{s\cdot \varphi \cdot r} = s\cdot \overline{\varphi}\cdot r$ for
$r,s\in R$ and $\varphi \in \mathit{Hom}_{\fieldk}(R,R)$.
Thus, $\varphi \in D^m_{\fieldk}(R)$ implies that
$\overline{\varphi} \in D^m_{\fieldk}(_RA)$.  Now use theorem \ref{T:DA=DRA}
to complete the lemma.
\end{proof}
\begin{remark}
By choice of $f_i$, we have
$\overline{\varphi} (r) = \varphi (r)$ for $r\in R$.  Hence the association
$\varphi \mapsto \overline{\varphi}$ is an injective map of $R$-bimodules.    
\end{remark}
\begin{theorem}\label{T:DA=ADRA}
$D^m_{\fieldk}(A)$ is generated as an $A$-bimodule by
$\{ \overline{\varphi} | \varphi \in D^m_{\fieldk}(R) \}$; that is,
\[D^m_{\fieldk}(A) = (A\tensor _R A^o)\cdot D^m_{\fieldk}(R) 
\cdot (A\tensor _R A^o).
\]
\end{theorem}
\begin{proof}
Let $\Phi \in D^m _{\fieldk}(A)$.  For each
$i,j \in \{0,1,2,\cdots , n \}$ let 
\begin{equation}\label{E:phiij}
\left( \Phi \right)_{i,j} = f_i \Phi \rho _{b_j}.
\end{equation}
Note that $\left( \Phi \right)_{i,j} (R) \subset R$.
For $r,s\in R$, we see that 
\[
r\cdot \left( \Phi \right)_{i,j} \cdot s
= \left( r\cdot \Phi \cdot s \right)_{i,j}.
\]
Hence, $\left( \Phi \right)_{i,j} \in D^m _{\fieldk}(R)$.
Now, by the dual basis lemma, 
\begin{align*}
\Phi (a) &= \sum _{j\geq 0} \Phi (f_j (a) b_j)\\
	&= \sum _{i,j \geq 0} f_i (\Phi (f_j(a)b_j))b_i.
\end{align*}
Hence, 
\begin{align}
\Phi &= \sum _{i,j \geq 0} \rho _{b_i} f_i \Phi \rho _{b_j} f_j \notag \\
     &= \sum _{i,j \geq 0} \rho _{b_i} \left( \Phi \right)_{i,j} f_j 
			\label{E:Phitophiij}
\end{align}
Since $f_j(A) \subset R$, we have
$\left( \Phi \right)_{i,j} f_j = \overline{\left( \Phi \right)_{i,j}} f_j$.
Hence, equation \ref{E:Phitophiij} gives
\[
\Phi = \sum _{i,j \geq 0} \rho _{b_i} \overline{\left( \Phi \right)_{i,j}}
			f_j
	\in (A\tensor _R A^o) \cdot \overline{D^m_{\fieldk}(R)} \cdot 
	(A\tensor _R A^o).
\]
Hence the theorem.
\end{proof}
\subsection{Ideals of $D_{\fieldk}(A)$}\label{S:ideals}
In this section, we show a one to one correspondence between ideals of
$D_{\fieldk}(A)$ and of $D_{\fieldk}(R)$.  
Let $\mathcal{I}_A$ and $\mathcal{I}_R$ denote the collection
of ideals in $A$ and $R$ respectively.
\begin{lemma}
For $I \in \mathcal{I}_A$, the set $f_0 I f_0$ is an ideal in 
$D_{\fieldk}(R)$.
\end{lemma}
\begin{proof}
The lemma follows from the fact that, for
$\varphi _1 , \varphi _2 \in D_{\fieldk}(R)$, and 
$\Phi \in D_{\fieldk}(A)$, we have $\varphi _1f_0 \Phi f_0 \varphi _2
= f_0 
\overline{\varphi _1} \Phi \overline{\varphi _2}f_0$.
\end{proof}
Define functions $\zeta$ and $\eta$ as follows:
\begin{align}
\zeta : \mathcal{I}_A &\to \mathcal{I}_R \notag \\
     I &\mapsto	f_0 I f_0 \label{E:I00}\\
\eta : \mathcal{I}_R &\to \mathcal{I}_A \notag \\
      J &\mapsto      D_{\fieldk}(A)\overline{J}D_{\fieldk}(A) \label{E:phibar}
\end{align}
\begin{theorem}
The correspondence $\zeta$ is a bijective function from $\mathcal{I}_A$ to
$\mathcal{I}_R$ with $\eta$ being its inverse function.
\end{theorem}
\begin{proof}
We show that $\eta \zeta$ is the identity on $\mathcal{I}_A$.
Clearly, $\eta (\zeta (I)) \subset I$.  Let any $\Phi \in I$. 
The $\left( \Phi \right) _{i,j}$ as defined in \ref{E:phiij} are in 
$I \cap f_0 I f_0$.
Referring to equation \ref{E:Phitophiij}, 
the claim is proved.

Now we show that $\zeta \eta$ is the identity on $\mathcal{I}_R$.
Again, it is obvious that $J \subset \zeta \eta (J)$.  For the reverse 
inclusion, use the fact that 
\[D_{\fieldk}(A) = (A\tensor _R A^o)
\overline{D_{\fieldk}(R)} (A\tensor _R A^o).
\]
For $a,b,c,d,p,q,m,n \in A, \psi _1, \psi _2 \in D_{\fieldk}(R)$ and
$\varphi \in J$,  we see that
\begin{align*}
f_0 &\left( (p \otimes q^o)  \overline{\psi _2} (c\otimes d^o)
	\overline{\varphi } (a \otimes b^o) \overline{\psi _1} 
 (m \otimes n^o) \right) f_0 \\
    &=
	\sum _{i,j,k \geq 0} 
	\left( f_0 (pb_k q) \cdot \psi _2 \cdot f_k (cb_i d)\right)
	\circ \left( \varphi \cdot f_i (a b_j b) \right)
	\circ \left( \psi_1 \cdot f_j(mn) \right) \\
     &\in J.
\end{align*}
Hence the theorem.
\end{proof}
\begin{corollary}
The ring $D_{\fieldk}(A)$ is noetherian if and only if the ring
$D_{\fieldk}(R)$ is.
\end{corollary}
\begin{corollary}
A ring $S$ is called Prime if for any ideals $P,Q$ of $S$, if $PQ = 0$ and 
$P \neq 0$, then $Q=0$.  The ring $D_{\fieldk}(A)$ is prime if and only
if $D_{\fieldk}(R)$ is.
\end{corollary}
\section{Differential operators on the Heisenberg algebras}\label{S:Heisen}
Let $\fieldk$ be a field and $n$ a positive integer.
Let $H_n$ denote the $n$th-Heisenberg algebra over $\fieldk$.  That is,
$H_n$ is a $\fieldk$-algebra with generators $h,x_1,x_2,\cdots ,x_n, y_1,
y_2,\cdots ,y_n$ such that
$[x_i , y_j] = \delta _{i,j} h$ and all the other commutators between the
generators equal 0.

In this section, we show that the ring of differential operators on
$H_n$ is generated by  two copies of $D_{\fieldk}(R)$ in the case of
zero characteristic and one copy in the non zero characteristic, where
$R$ denotes the centre of $H_n$.   Note that in the case of non zero 
characteristic, the centre is very large (that is, $H_n$ is free of finite
rank as a module over its centre).

\subsection{Characteristic of \fieldk is 0}
In this case, the centre of $H_n$ is $\fieldk [h]$, the polynomial ring
in one variable. Here, we have two different inclusions of
$D_{\fieldk}(\fieldk [h])$ into $D_{\fieldk}(H_n)$. 

Let $I = (i_1, i_2, \cdots , i_n) \in (\mathbb{Z}_+ )^n$ be a multi index.
Denote by $\mathbf{x} ^I$ the element $x_1 ^{i_1} x_2^{i_2} \cdots x_n ^{i_n}$.
Note that every element $a \in H_n$ can be written uniquely as
$a = \sum _{I,J} p_{I,J} (h) \mathbf{x}^I \mathbf{y}^J$ where
$I,J$ are multi indices in $(\mathbb{Z}_+)^n$ and $p_{I,J}(h)$ is a 
polynomial in $h$ with coefficients in $\fieldk$.

For a multiindex $I = (i_1,i_2, \cdots , i_n)$, let
$|I|$ denote the sum $(i_1 + i_2 +\cdots + i_n)$. We define two kinds of
degree on $H_n$.\hfill
\begin{enumerate}
\item	For $a \in H_n$  such that $a = \sum p_{I,J} \mathbf{x}^I \mathbf{y}^J$
	define $\mathit{deg}_1$ of $a$ as
	\[
	\mathit{deg}_1 (a) = \textit{max} \{ |I| + |J| \mid p_{I,J} \neq 0 \}.
	\]
\item	Define $\mathit{deg}_2 (x_i) = \mathit{deg}_2 (y_i ) = 1$
	and $\mathit{deg}_2(h) = 2$ and extend this degree to the entire ring.
\end{enumerate}

We see that $H_n$ is filtered as $H_n = \cup _{k \geq 0} H_n ^k$
where 
\[
	H_n ^k = \{ a | \mathit{deg}_1 (a) \leq k \}.
\]
\begin{note}\label{N:Hn}
For $a \in H_n^k$, and $r \in \{ x_1,x_2,\cdots , x_n, y_1,y_2,\cdots
, y_m \}$, we have
$[r,a] \in H_n^{k-1}$.
\end{note}
\begin{lemma}
For any $H_n$-bimodule $M$, let $m \in M_{\mathit{diff}}$ (as defined in
definition \ref{D:diff}). Then there exists a $k \geq 0$ \\
such that
$[ \cdots [[m, r_1],r_2],\cdots ,r_k] =0$ for\\
$r_i \in \{x_1,x_2,\cdots ,x_n, y_1,y_2,\cdots ,y_n \}$.
\end{lemma}
\begin{proof}
Let $m \in \mathcal{Z}_t M$ (definition \ref{D:diff}) for some $t \geq 0$ 
such that $ m = a.n$ for some $a \in H_n$ and $ n \in \mathcal{Z}(M \diagup
\mathcal{Z}_{t-1}M)$.  It is enough to show that there exists an $l \geq 0$ 
such that $[\cdots [[m,r_0],r_1],\cdots ,r_l] \in \mathcal{Z}_{t-1}M$
for \\$r_i \in \{x_1,x_2,\cdots ,x_n, y_1,y_2,\cdots ,y_n \}$.
If $a \in \fieldk [h]$ then $l=1$.  Else, $ a \in H_n^l$ for some $l\geq 0$.
By referring to the note \ref{N:Hn}, we have the lemma.
\end{proof}
\begin{corollary}
Let $\varphi \in D_{\fieldk}(H_n)$.  Then there exists a $k\geq 0$ such that
$[\cdots [[\varphi , r_0],r_1],\cdots ,r_k] = 0$ for \\
$r_i \in \{x_1,x_2,\cdots ,x_n, y_1,y_2,\cdots ,y_n \}$.
Note that a $\varphi \in \mathit{Hom}_{\fieldk}(H_n,H_n)$ satisfying
this condition is in $D_{\fieldk}(H_n)$.
\end{corollary}
The corollary above provides another filtration of $D_{\fieldk}(H_n)$ 
given by $D_{\fieldk}(H_n) = \cup _{l\geq 0}M_l$ where
\[
M_l = \{ \varphi \in D_{\fieldk}(H_n) |
		[\cdots [[ \varphi ,r_0],r_1],\cdots ,r_l ]=0 \}
\]
for
$r_i \in \{x_1,x_2,\cdots ,x_n,y_1,y_2,\cdots ,y_n \}$.
\begin{note}
$M_l$ is closed under $+$ and $[M_l ,h] \subset M_{l-2}$.
\end{note}
\begin{lemma}
$M_l \cdot M_s \subseteqq M_{l+s}$
\end{lemma}
\begin{proof}
Immediate once we see that $[\varphi _1 \varphi _2 , r]
= \varphi _1 [\varphi _2 ,r] + [\varphi _1 ,r]\varphi _2$.
\end{proof}
\begin{definition}\label{D:diffs}
Let $\partial _{x_l} ,\partial _{y_l},\partial _h ,\overline{\partial _h}\in 
\mathit{Hom}_{\fieldk}(H_n,H_n)$ be defined as
\begin{align*}
\partial _{x_l} (p x_1^{i_1} x_2 ^{i_2} \cdots x_n^{i_n}
		y_1^{j_1}y_2^{j_2} \cdots y_n ^{j_n})
		&= i_l p x_1^{i_1} x_2^{i_2}\cdots x_l^{i_l -1} \cdots
		x_n^{i_n} y_1^{j_1}y_2^{j_2} \cdots y_n ^{j_n},\\
\partial _{y_l} (p x_1^{i_1} x_2 ^{i_2} \cdots x_n^{i_n}
		y_1^{j_1}y_2^{j_2} \cdots y_n ^{j_n})
		&= j_l p x_1^{i_1} x_2^{i_2}\cdots 
		x_n^{i_n} y_1^{j_1}y_2^{j_2} \cdots y_l^{j_l -1} \cdots
		y_n ^{j_n},\\
\partial _h     (p x_1^{i_1} x_2 ^{i_2} \cdots x_n^{i_n}
		y_1^{j_1}y_2^{j_2} \cdots y_n ^{j_n})
		&= p^{\prime} x_1^{i_1} x_2 ^{i_2} \cdots x_n^{i_n}
		y_1^{j_1}y_2^{j_2} \cdots y_n ^{j_n},\\
\overline{\partial _h}     (p y_1^{i_1} y_2 ^{i_2} \cdots y_n^{i_n}
		x_1^{j_1}x_2^{j_2} \cdots x_n ^{j_n})
		&= p^{\prime} y_1^{i_1} y_2 ^{i_2} \cdots y_n^{i_n}
		x_1^{j_1}x_2^{j_2} \cdots x_n ^{j_n},
\end{align*}
where $p^{\prime }$ denotes the usual derivative of $p$ with respect to $h$.
\end{definition}		
We list some immediate properties:\hfill
\begin{enumerate}
\item	$[\partial _r, \partial _s] = 0$ for $r,s \in \{x_1, \cdots 
	,x_n,h,y_1,\cdots , y_n \}$.
\item 	$[\partial _{x_l} , x_l] = 1$ and $[\partial _{x_l},r] =0$ for
	\[
	r \in \{x_1,\cdots ,x_n,h,y_1,\cdots ,y_n \} \setminus \{x_l \}.
	\]
\item	$[\partial _{y_l} , y_l] = 1$ and $[\partial _{y_l},r] =0$ for
	\[
	r \in \{x_1,\cdots ,x_n,h,y_1,\cdots ,y_n \} \setminus \{y_l \}.
	\]
\item 	$[\partial _h , h] = 1$ , $[\partial _h, y_l] = - \partial_{x_l}$
	and $[\partial _h ,x_l] = 0$, for $1\leq l \leq n$.
\item	The above properties show that $\partial _{x_l}, \partial _{y_l}
	\in M_1$ and $\partial _h \in M_2$.
\item	$\lambda _{x_l} - \rho _{x_l} = h \partial _{y_l}$, and 
	$\rho _{y_l} - \lambda _{y_l} = h \partial _{x_l}$.
\item	$\overline{\partial _h} = \partial _h + \sum_l
	\partial _{x_l} \partial _{y_l}$.
\item	$[\overline{\partial _h},x_l]= \partial _{y_l}$, 
	$[\overline{\partial _h},y_l]= 0$ and
	$[\overline{\partial _h},h]= 1$.
\end{enumerate}
Following the theorem 2.3.4 of \cite{LR}, we show the following
\begin{proposition}
Let characteristic of \fieldk be 0.
The \fieldk -algebra \\
$D_{\fieldk}(H_n)$ is generated by left multiplications
by elements of $H_n$ and 
by 
\[
\{ \partial _{x_l} , \partial _{y_l} \}_{1\leq l \leq n} , \partial _h.
\]
\end{proposition}
\begin{proof}
Let $R_s \subset H_n$ denote the \fieldk -span of monomials in 
\[
x_1,\cdots ,x_n,h,y_1,\cdots ,y_n \textit{ of }\mathit{deg}_2 \leq s.
\]
\textit{Claim :} Let $D\in M_s$ be such that $D |_{R_s} =0$.
Then $D= 0$.\\
When $D \in M_0$ then $D = \rho _{D(1)}$ and hence the claim. Assume that
we have proved the claim for $s \lneq i$ and fix $D \in M_i$ such that
$D|_{R_j} =0$ for some $j\geq i$.
It is enough to show that for $c\in R_j, c^{\prime} \in R_{j-1}$, we have
$D(x_lc) = D(y_l c) = D(h c^{\prime}) = 0$. Note that
$D(x_lc) = [D,x_l] (c) + x_l D(c) = 0+0$.
Similar argument for the variables $y_l$ and the fact that
$[M_i ,h] \subset M_{i-2}$ complete the claim.

Let $A\subset D_{\fieldk}(H_n)$ be the \fieldk -subalgebra generated
by $H_n$ and \\
$\{ \partial _{x_l} , \partial _{y_l} \}_{1\leq l \leq n} , 
\partial _h$.\\
\textit{Claim :} $H_n$ is a simple $A$-module.\\
By assumption, the characteristic of the field is 0.   Hence the claim follows
from the fact that given a $0 \neq c \in R_j$, there exists a 
$D \in \{ \partial _{x_l} , \partial _{y_l} \}_{1\leq l \leq n} \cup 
\{ \partial _h \}$ such that $D(c) \neq 0$ and $D(c) \in R_{j-1}$.

Note that $\mathit{Hom}_A(H_n,H_n) = \fieldk$.  Now fix $ 0 \neq D \in M_s$.
Then by the Jacobson Density Theorem, we can find $d \in A$ such that
$d |_{R_s} = D|_{R_s}$.  If $d \notin M_s$, then clearly, $d|{R_s} = 0$.
Now, by the first claim, $d =D$.  Hence the proposition.
\end{proof}
\begin{corollary}\label{C:weyl}
If $A_n$ denotes the $n$th-Weyl algebra (that is $h =1$)
over a field of characteristic
0, then $D_{\fieldk}(A_n) = A_{2n}$ 
\end{corollary}
\begin{proof}
By the above proposition, $D_{\fieldk}(A_n)$ is generated by 
$\{ \partial _{x_l},\partial _{y_l} \}_{1\leq l \leq n}$ and left
multiplication by elements of $A_n$.  Since
$[x_i,y_j] = \delta _{i,j}$, we have $\partial _{x_l} , \partial _{y_l}$
are inner (A more general statement is true, due to Dixmier (Lemma 4.6.8 
of \cite{D}):  
All 
derivations on a Weyl algebra are inner).
That is, $D_{\fieldk}(A_n) = D_{\fieldk}^0 (A_n)$.
By corollary \ref{C:D0}, we have
a surjection  $A_n \tensor A_n^o \to D^0_{\fieldk}(A_n)$.
Note that $A_n^o $ is isomorphic to $A_n$ by mapping
$x_l^o \mapsto  -y_l$ and $y^o_l\mapsto x_l$.
Also, $A_n \tensor A_n \cong A_{2n}$ by Corollary 1.2, page 122 of \cite{C}.
Thus, we have a surjection $A_{2n} \to D_{\fieldk}^0(A_n)$.  Now use the fact
that $A_{2n}$ is simple to complete the corollary.
\end{proof}
\begin{theorem}\label{T:Hnchar0}
Let characteristic of \fieldk be 0.
The ring $D_{\fieldk}(H_n)$ is generated by left multiplication by elements
of $H_n$ and by $\{\partial _h , \overline{\partial _h} \}$.
That is, the ring $D_{\fieldk}(H_n)$ is generated by two copies of
$D_{\fieldk}(\fieldk [h])$ and inner derivations.
\end{theorem}
\begin{proof}
By the properties following definitions \ref{D:diffs} we see that
$\partial _h$ and $\overline{\partial _h}$ generate $\partial _{x_l}$ and
$\partial _{y_l}$ for all $l$.
The previous proposition completes the theorem.
\end{proof}
\begin{theorem}\label{T:DHnsimple}
Let \fieldk be a field of characteristic 0.  The ring of differential 
operators on $H_n$ is simple.
\end{theorem}
\begin{proof}
Let $a \in D_{\fieldk}(H_n)$.  Then $a$ can be written as
a \fieldk -linear combination of monomials of the form
$h^m \mathbf{x}^I \mathbf{y}^J \partial ^s_h \mathbf{\partial _x}^K 
\mathbf{\partial _y}^L$
where 
\[
\mathbf{\partial_x}^K = \partial _{x_1}^{k_1} \partial _{x_2}^{k_2}
\cdots \partial _{x_n}^{k_n}
\]  
where $K = (k_1,k_2,\cdots ,k_n)$ a multiindex. 
Let $\mathcal{I}$ be an ideal in $D_{\fieldk}(H_n)$.  
Let $0 \neq a \in \mathcal{I}$.  As $[\partial _h ^s ,h] = s \partial _h^{s-1}$
 and the fact that $h$ commutes with all the other generators,
we can assume that $\partial _h$ does not appear in the expression of
$a$. Now use the fact that 
$[x_l^k \partial _{y_l}^s , y_l] = khx_l^{k-1} \partial _{y_l}^s
				+ sx_l^k\partial _{y_l}^{s-1}$ 
and the fact that $y_l$ commutes with all the other generators, 
to assume that in the expression of $a$, the generators $x_l$ and 
$\partial _{y_l}$ do not appear.  Similarly, as
$[y_l^k \partial _{x_l}^s ,x_l] = -khy_l^{k-1} \partial _{x_l}^s
				+ sy_l^k\partial _{x_l}^{s-1}$, we can
assume that $a$ is a polynomial in $h$.  Now use the fact that 
$[h^s,\partial_h] = sh^{s-1}$ to conclude that there is a non zero scalar in
$\mathcal{I}$ and hence $\mathcal{I} =D_{\fieldk}(H_n)$.
\end{proof}
\subsection{Characteristic of \fieldk = $p\neq 0$.}
In this case, the centre is the polynomial ring in $2n+1$ variables
\[
R := \fieldk [ h, x_1^p, x_2^p,\cdots ,x_n^p, y_1^p,y_2^p,\cdots ,y_n^p].
\]

Theorem 2 of \cite{R} shows that the $n$th-Weyl algebra is Azumaya over
its centre when characteristic of \fieldk is nonzero.  The same proof works
to show that $H_n$ is Azumaya over its centre.  Now we refer to the section
on Azumaya algebra to claim: 
\begin{theorem}
$D_{\fieldk}(H_n) = (H_n \tensor _R H_n^o) \cdot D_{\fieldk}(R) \cdot
			(H_n \tensor _R H_n^o)$.
\end{theorem}
In \cite{S}, the differential operators on polynomial ring in one variable, 
on a field of nonzero characteristic has been studied.  In particular
it is shown that $D_{\fieldk}(R)$ is simple.
\begin{corollary}\label{C:Hncharp}
Let \fieldk be a field of non zero characteristic.  The ring 
$D_{\fieldk}(H_n)$ is simple.
\end{corollary}

\section{Concluding remarks and acknowledgements}\label{S:CR}
This work suggests that if $R$ is the centre of $A$, and if 
there is a way of embedding $D_{\fieldk }(R)$ into 
$D_{\fieldk}(A)$, then $D_{\fieldk}(R)$ generates $D_{\fieldk}(A)$ as an 
$A^e$-module. Further natural questions are to find differential operators on 
the enveloping algebras of Lie algebras.

This work was part of my thesis written at
Indiana University, Bloomington, Indiana, U.S.A, under the guidance
of Professor Valery A. Lunts.  I would like to thank
Professors Darrell Haile and
Valery Lunts for their generous help and suggestions.  I would also like to
thank Professors Dipendra Prasad and R.Sridharan for suggesting some useful
questions.  I would like to thank Dr. Timothy McCune for useful discussions.

 \end{document}